\documentclass[10pt]{amsart}
\usepackage{amsmath,amsfonts,amscd,amssymb}

\title{Quotients of functors of Artin rings}
\author{Tim Dokchitser}
\date{November 21, 2005}
\address{Tim Dokchitser\vskip 0mm
Robinson College\vskip 0mm
Cambridge, CB3 9AN\vskip 0mm
United Kingdom}
\email{t.dokchitser@dpmms.cam.ac.uk}

\theoremstyle{plain}
\newcounter{thmcount}
\newtheorem{theorem}[thmcount]{Theorem}

\newtheorem{lemma}[thmcount]{Lemma}
\theoremstyle{definition}

\newtheorem*{remark*}{Remark}

\newtheorem*{question*}{Question}

\def\newmathop#1{\expandafter\gdef\csname #1\endcsname{\mathop{\rm #1}\nolimits}}

\newmathop{id}
\newmathop{Im}
\newmathop{Hom}
\newmathop{Aut}

\def\surjects{\rightarrow\penalty 1000\kern-0.8em\penalty 1000\rightarrow}
\def\leftsurjects{\leftarrow\penalty 1000\kern-0.8em\penalty 1000\leftarrow}
\def\compose{{\raise 1pt\hbox{$\scriptscriptstyle\circ$}}}

\def\injects{
   \mbox{
     $\longrightarrow$
     \kern-\mathsurround
     \kern-\mathsurround
     \kern-2.1em\raise0.29em\hbox{{$\scriptstyle\subset$}}
     \kern1.1em
   }
}

\def\overarrow#1{{\buildrel #1 \over \longrightarrow}}

\let\W\Lambda

\def\iW{_\Lambda}
\def\iA{_A}
\def\iB{_B}
\def\iF{_F}

\def\iC{_C}

\def\cD{{\mathcal D}}
\def\cF{{\mathcal F}}

\def\cG{{\mathcal G}}
\def\si#1{_{#1}}

\def\beq{\begin{eqnarray}}
\def\eeq{\end{eqnarray}}
\def\ben{\begin{enumerate}}
\def\een{\end{enumerate}}
\def\bit{\begin{itemize}}
\def\eit{\end{itemize}}

\def\Art{{Art}\iW}                             % Artinian augm. Lambda-algebras
\def\Sets{{S\hskip -0.1em et\hskip -0.03em s}}

\let\lar\longrightarrow
\let\iso\cong

\def\hF{h\iF}

\def\Auttau{\Aut^T\iW}
\newmathop{id}

\begin{document}

\rlap{\hskip 30cm.}\vskip -0.3cm

\maketitle

%arXiv Abstract: 
%In infinitesimal deformation theory, a classical criterion due to %Schlessinger gives an intrinsic characterisation of functors that are %pro-representable, and more generally, of the ones that have a hull. Our %result is that in this setting the question of characterising group %quotients can also be answered. In other words, for functors of Artin rings %that have a hull, those that are quotients of pro-representable ones by a %constant group action can be described intrinsically.

One of the fundamental problems in the study of moduli spaces is to
give an intrinsic characterisation of representable functors of schemes,
or of functors that are quotients of representable ones of some sort.
Such questions are in general hard,
leading naturally to geometry of algebraic stacks and spaces (see \cite{Art,KM}).

On the other hand, in infinitesimal deformation theory a classical
criterion due to Schlessinger \cite{Schl}
does describe the pro-representable functors and, more generally,
functors that have a hull.
Our result is that in this setting the question of describing
group quotients also has a simple answer. In other words, for functors of Artin rings that
have a hull, those that are quotients of
pro-representable ones by a constant group action can be described
intrinsically.

%Introduction: representable functors / stacks; quotients.
%
%Deformation functors are quotients.
%
%Conjecture: everything is a quotient.
%
%Application: tangent space
%
%Two cases: constant group functor, pro-representable formally smooth
%group functor.
%
%Result: necessary \& sufficient conditions for quotients by groups.
%Plus conjecture for the other case.

To set up the notation, let $\W$ be a complete Noetherian local ring with residue field $k$, and
fix an isomorphism $\W/m\iW\iso k$.
We write $\Art$ for the category of Artinian local $\W$-algebras $A$
given together with an augmentation $A/m\iA\iso k$.
Morphisms (denoted $\Hom\iW(A,B)$) are local homomorphisms of $\W$-algebras that commute
with the augmentations.

Every complete Noetherian local $\W$-algebra $F$ with an augmentation
(not necessarily Artinian) gives rise to a covariant functor
$$
  \hF = \Hom\iW(F,-): \>\Art\lar\Sets,
$$
and a functor $\cF$ isomorphic to some $\hF$ is said to be pro-representable.
A functor $\cG$ has a hull if there is a formally smooth natural
transformation $\alpha: \hF\to\cG$ for some~$F$ such that $\alpha$ is bijective
on rings of the form $k[V]$ with $V$ a $k$-vector space.

Recall that by Schlessinger's theorem (\cite{Schl}, Thm 2.11),
$\cF$ is pro-representable if and only if $\cF(k)$ consists of one element,
the fibre product map
\beq\label{e:fib}
  \cF(A\times\iB C) \>\>\lar\>\> \cF(A)\times\si{\cF(B)}\cF(C)
\eeq
is bijective for all $A\to B\leftsurjects C$ in $\Art$,
and the tangent space $T\cF=\cF(k[t]/t^2)$ is finite-dimensional.
Also, $\cF$ has a hull if and only if the same conditions hold except that
\eqref{e:fib} is only assumed to be surjective.

For a complete Noetherian local $\W$-algebra $F$,
an action of a group $\Gamma$ on $\hF$ is given by a homomorphism
$\Gamma\to\Aut\iW(F)$, and such an action naturally gives rise to
a quotient functor $\cF/\Gamma$. We write
$TF=\Hom_k(m\iF/(m\iF^2+m\iW F),k)$ for the tangent space of $F$ over $W$
(so $TF=T\hF$ naturally) and set
$$
  \Auttau(F)=\ker\bigl(\Aut\iW(F)\to\Aut\iW(TF)\bigr).
$$

\begin{theorem}\label{t:qgmain}
For a functor $\cD:\Art\to\Sets$
the two conditions are equivalent:
\bit
  \item[(i)] There exists a pro-representable functor $\cF\iso\Hom\iW(F,-)$ and a subgroup
  $\Gamma\subset\Auttau(F)$ such that $\cD\iso\cF/\Gamma$.
  \item[(ii)] $\cD$ possesses a hull and $\cD(A\times\iB C)\surjects \cD(A)\times_{\cD(B)}\cD(C)$
    for all $A\rightarrow B\leftarrow C$.
%  In this case $\cF$ is a hull of $\cD$ and $\Gamma$ is unique
%  up to conjugacy. %namely $\Gamma=\Aut\si{\cD}(\cF)\subset\Auttau(F)$.
\eit
\end{theorem}
\begin{proof}
The following lemma lists some properties of quotients of pro-representable
functors by group actions and, in particular, proves the easy implication,
(i) $\Rightarrow$ (ii).

\begin{lemma}\label{t:qg}
Let $\cF=h\iF$ and let $\cD=\cF/\Gamma$ with $\Gamma\subset\Aut\iW(F)$. Then
\ben
  \item The quotient map $q: \cF\to\cD$ is formally smooth.
  \item The map $\cD(A\times\iB C)\to\cD(A)\times_{\cD(B)}\cD(C)$
    is surjective for all $A\rightarrow B\leftarrow C$.
  \item $\cD$ has a hull if and only if $\Gamma\subset\Auttau(F)$,
  in other words, if $\Gamma$ acts trivially on the tangent space of $F$.
  \item $\cD$ is pro-representable if and only if $\Gamma=\{1\}$.
\een
\end{lemma}

\begin{proof} (1) For $\pi: A\surjects B$ in $\Art$, we need to show that
$$
  \cF(A) \surjects \cF(B)\times_{\cD(B)}\cD(A) \>.
$$
Let $a\in\cD(A)$ and $\tilde b\in\cF(B)$ be such that $\pi(a)=q(\tilde b)$
in $\cD(B)$, and lift $a$ to $\tilde a\in\cF(A)$.
Then $g\cdot \pi(\tilde a)=\tilde b$ for some $g\in \Gamma$, and
$g\cdot\tilde a$ has the required properties.

\noindent (2)
Let $\pi: A\rightarrow B$ and $\rho: C\rightarrow B$ in $\Art$. Let
$a \in\cD(A)$ and $c\in\cD(C)$ be such that
$\pi(a)=\rho(c)$ in $\cD(B)$ and lift them to some $\tilde a\in\cF(A)$ and
$\tilde c\in\cF(C)$. Then
$\pi(\tilde a)$ and $\rho(\tilde c)$ in $\cF(B)$ map to the same element
in $\cD(B)$, hence there is a $g\in \Gamma$ such that
$$
  g\cdot\pi(\tilde a) = \rho(\tilde c) \>.
$$
Replace $\tilde a$ by $g\cdot\tilde a$, so that $\tilde a$ still maps to
$a\in\cD(A)$, but now $\pi(\tilde a)=\rho(\tilde c)$.
Since $\cF$ commutes with fibre products, there is
$\tilde r\in\cF(A\times\iB C)$
projecting to $\tilde a\in\cF(A)$ and $\tilde c\in\cF(C)$.
Then its image $r\in\cD(A\times\iB C)$ is a required lift of $(a,c)$.

\noindent(3)
The ``if'' part follows from the fact that $\cF$ has a tangent space
and that $\cF(k[V])\to\cD(k[V])$ is bijective for all $k$-vector spaces $V$
by the assumption on~$\Gamma$.

For the converse, assume that the action of $\Gamma$ on $TF$ is non-trivial
but the quotient $\cD=\cF/\Gamma$ has a hull.
Write $V=\cF(k[t]/t^2)=\cF(k[\epsilon])$. Let
$k[\epsilon_1,\epsilon_2]$ denote the ring $k[t_1,t_2]/(t_1^2,t_2^2,t_1t_2)$
and consider the map
$$
  \pi: \cD(k[\epsilon_1,\epsilon_2]) \lar \cD(k[\epsilon]) \times \cD(k[\epsilon]),
$$
whose components $\pi_1$ and $\pi_2$ are induced by the natural projections.
Since $\cD$ has a hull, it has a tangent space,
so $\pi$ is a bijection. Now,
$$
  \cD(k[\epsilon_1,\epsilon_2]) = \cF(k[\epsilon_1,\epsilon_2])/\Gamma = (V\oplus V)/\Gamma \>,
$$
and
$$
  \cD(k[\epsilon]) \times \cD(k[\epsilon]) = (V/\Gamma) \times (V/\Gamma) \>.
$$
Moreover, the action of $\Gamma$ on $\cF(k[\epsilon_1,\epsilon_2])=V\oplus V\>$
is diagonal,
$$
  g \cdot (v_1,v_2) = (g\cdot v_1,g\cdot v_2), \qquad v_1,v_2\in V\>, g\in \Gamma\>,
$$
by compatibility of the action with the two inclusions
$k[\epsilon]\injects k[\epsilon_1,\epsilon_2]$.
Since the action of $\Gamma$ on $V$ is non-trivial, there are
$v_1\ne v_2\in V$ such that $g\cdot v_1=v_2$ for some $g\in \Gamma$.
Then
$$
  h\cdot (v_1,v_1) = (h\cdot v_1,h\cdot v_1) \ne (v_1,v_2)
$$
for any $h\in\Gamma$.
Thus $(v_1,v_1)$ and $(v_1,v_2)$ give two distinct
elements of $\cD(k[\epsilon_1,\epsilon_2])$, but
$\pi_1(v_1,v_2)=\pi_2(v_1,v_2)$ in $\cD(k[\epsilon]) \times \cD(k[\epsilon])$,
contradicting injectivity of $\pi$.

\noindent(4)
If $\Gamma=\{1\}$, then $\cD=\cF$ is pro-representable. Conversely,
if $\cD$ is pro-representable, it certainly has a hull, so
$\Gamma\subset\ker\bigl(\Aut\iW(F)\to\Aut(TF)\bigr)$ by part (3).
So $\cF\to\cD$ is a formally smooth map of pro-representable functors,
which is an isomorphism on the tangent spaces, so it is an isomorphism.
\end{proof}

As for the implication (ii) $\Rightarrow$ (i) of the theorem,
to show that $\cD\iso\cF/\Gamma$ for some $\Gamma$ we need
to exhibit sufficiently many automorphisms commuting with the projection $\cF\to\cD$.
This has two ingredients, an ``automorphism lemma'' below, and the fact
that $\cD$ takes injections to injections which is a consequence
of the assumed surjectivity $\cD(A\times\iB C)\surjects \cD(A)\times_{\cD(B)}\cD(C)$
for all $A\to B\leftarrow C$.

\begin{lemma}
\label{t:autlemma}
Let $q: \cF\to\cD$ be formally smooth with
$\cF=\Hom\iW(F,-)$ pro-representable.
Take $x,y\in\cF(A)$ with $q(x)=q(y)$ and such that
$y$ is surjective, when considered as a homomorphism $F\to A$. Then
there exists a natural transformation $\alpha:\cF\to\cF$ for which
$$
  \begin{array}{ccc}
   \cF\!\!\!\!\!& {\buildrel\alpha\over{\hbox to 35pt{\rightarrowfill}}}&\!\!\!\cF \cr
                &\llap{$\scriptstyle q\!$}\searrow
                 \qquad\swarrow\rlap{$\scriptstyle\!q$}&\cr
                & \!\!\cD\!\!                                  &\cr
  \end{array}
$$
commutes and such that $\alpha(x)=y$.
\end{lemma}

\begin{proof}
Consider a commutative diagram,
$$
\begin{array}{rcl}
                  \cF(F) & \overarrow{q} & \cD(F)\ni q           \cr
   \raise 2pt\hbox{$\scriptstyle\cF(x)$}
   \Bigl\downarrow^{\vphantom{Z}}_{\vphantom{y}}\quad
                         &                    & \quad\Bigl\downarrow\rlap{\raise 2pt\hbox{$\scriptstyle
                                        \scriptstyle\cF(x)$}}\cr
              y\in\cF(A) & \overarrow{q} & \cD(A)                   \cr
\end{array} \>.
$$
viewing $q$ both as a natural transformation and as an element of $\cD(F)$.
As $x$ is surjective and $\cF\to\cD$ is formally smooth, there exists
$g\in\cF(F)$ lifting $(q,y)$.

It is now a tautology that $g\in\Hom\iW(F,F)$ is a homomorphism,
whose associated natural transformation $\alpha$ has the required properties.
First,
$q(g)=q$ implies that $q\alpha=q$ as natural transformations.
Second, $\cF(x)(g)=y$ says precisely that $\alpha(x)=y$.
\end{proof}

We can now complete the proof of the theorem.
Let $q: \hF=\cF\to\cD$ be a hull of $\cD$ and let $\Gamma\subset\Aut\iW(F)$
consist of those automorphisms $g$ which, considered as elements of $\Aut(\cF)$,
satisfy $q g=q$ (as natural transformations). Since $q$ is an isomorphism
on tangent spaces, $\Gamma\subset\ker\bigl(\Aut\iW(F)\to\Aut(TF)\bigr)$,
and it suffices to prove that $\cD\iso\cF/\Gamma$. Clearly $q$
factors through $\cF/\Gamma$ and
$$
  \cF(A)/\Gamma \to \cD(A)
$$
is surjective for all $A\in\Art$, since $\cF(A)\surjects\cD(A)$ by
formal smoothness.
To prove injectivity, assume $x,y\in\cF(A)$ are such that
$q(x)=q(y)\in\cD(A)$. We claim that $g\cdot x=y$ for some $g\in \Gamma$.

Consider $x$ and $y$ as homomorphisms $F\to A$. We first
reduce to the case that $x,y$ are surjective. Let $A'\subset A$ be the
$\W$-subalgebra generated by $\Im x$ and $\Im y$. Then both $x$
and $y$ factor through $A'$,
$$
  x,y: F \lar A' \injects A \>.
$$
In other words $x,y\in\cF(A)$ lie in the image of $\cF(A')\injects\cF(A)$.
Let $x',y'\in\cF(A')$ be the same homomorphisms, considered as elements
of $\cF(A')$. We claim that $q(x')=q(y')\in\cD(A')$.

We know that $q(x')$ and $q(y')$ have the same image in $\cD(A)$.
By the second assumption on $\cD$, the map
$$
  \cD(A') = \cD(A'\times\iA A') \to \cD(A') \times_{\cD(A)}\cD(A')
$$
is surjective. Equivalently  $\cD(A')\injects\cD(A)$, so
$\cD$ takes injections to injections, implying $q(x')=q(y')$. If we can
find $g\in \Gamma$ for which $g\cdot x'=y'$, then $g\cdot x=y$ as required.
So we can replace $A$ by $A'$, in other words assume that $A$ is generated
by $\Im x$ and $\Im y$ as a $\W$-algebra.

We claim that in this case both $x$ and $y$ are surjective.

Indeed, let $B=\Im(x)$ and $C=\Im(y)$.
As $A$ is generated by $B$ and $C$ as a $\W$-algebra,
the cotangent space $V=m\iA/(m\iA^2+m\iW A)$ is generated, as a vector space,
by $m\iB V$ and $m\iC V$. Thus, if we show that $m\iB V = m\iC V$, then
it follows that $x,y$ are surjective on cotangent spaces, hence surjective.
%?reference
%(Remark \ref{r:nak}).

Consider the projection $A\to A/m\iA^2$, composed with $x$ and $y$:
$$
  F \>\>\overarrow{x,y}\>\> A \>\>\surjects\>\> A/m\iA^2 \>.
$$
The compositions $\bar x$ and $\bar y$ define elements
$\bar x,\bar y\in\cF(A/m\iA^2)$. Since $q(\bar x)=q(\bar y)$,
and $q: \cF\to\cD$ is a bijection on the rings of the form $k[V]$
for $V$ a $k$-vector space (such as $A/m\iA^2$), it follows that $\bar x=\bar y$. So $m\iB V=m\iC V$
and $x,y$ are both surjective.

%In summary, we have surjections $x,y:F\to A$ and $q(x)=q(y)\in\cD(A)$
%if $x,y$ are considered as elements of $\cF(A)$. We will now exhibit
%$g\in \Gamma$ with $g\cdot x=y$.

Now, by Lemma \ref{t:autlemma}
there exists a homomorphism $g: F\to F$
such that the corresponding natural transformation $\alpha:\cF\to\cF$
commutes with $q$ and such that $xg=y$. Since $q\alpha=q$, it follows
that $g$ is identity on the tangent space of $F$. In particular, it is
an automorphism of $F$ and $g\in \Gamma$. Also $g\cdot x=y$, as asserted.
This completes the proof of Theorem \ref{t:qgmain}.
\end{proof}

\begin{remark*}
Rather than considering a group action, take more generally a formally
smooth group functor
$\cG:\Art\to\Sets$ with an action $\cG\times\cF\to\cF$ on a
pro-representable functor $\cF$. If the induced action on the tangent space
of $\cF$ is trivial, it is not hard to see that the quotient $\cF/\cG$ by
this action has $\cF$ as a hull, as in the case when $\cG$ is a ``constant
group functor'' as above. It would be interesting to know whether,
conversely, every $\cD$ that has a hull can be represented by such a
quotient and, if not, whether such $\cD$ can be characterised intrinsically.
Thus we end with the following question:
\end{remark*}

\begin{question*}
Assume that $\cD: \Art\to\Sets$ has a hull $q: \hF\to\cD$. Does there exist
a formally smooth group functor $\cG$ acting on $\hF$ such that
$\cF/\cG\iso\cD$ with $q$ as the quotient map?
\end{question*}

\end{document}